\documentclass[12pt]{amsart}
\usepackage{amsmath,amsxtra,amssymb,amsfonts}

\numberwithin{equation}{section}

\newcommand{\mel}{\end{eqnarray*}}
\def\fr{\begin{align*}}

\newcommand{\kl}{\pl \le \pl}

\newcommand{\lel}{\pl = \pl}

\newcommand{\nz}{{\mathbb N}}

\newcommand{\rz}{{\mathbb R}}

\newcommand{\cz}{{\mathbb C}}

\newcommand{\ten}{\otimes}

\newcommand{\p}{\varphi}
\newcommand{\pl}{\hspace{.1cm}}

\renewcommand{\a}{\alpha}
\newcommand{\al}{\alpha}
\newcommand{\si}{\sigma}

\newcommand{\eps}{\varepsilon}

\newcommand{\M}{{\mathcal M}}

\renewcommand{\P}{\mathcal P}

\newcommand{\N}{{\mathcal N}}

\newcommand{\noo}{\left \|}
\newcommand{\rrm}{\right \|}

\newcommand{\intt}{\int\limits}

\newcommand{\NN}{{\mathcal N}}

\newcommand{\OL}{\mathcal {OL}}
\newcommand{\RR}{{\mathcal M}}

\newtheorem{theo}{Theorem}

\newtheorem{lemma}{Lemma}[section]
\newtheorem{prop}[lemma]{Proposition}
\newtheorem{theorem}[lemma]{Theorem}

\newtheorem{cor}[lemma]{Corollary}

\newtheorem{rem}[lemma]{Remark}

\newcommand{\re}{\begin{rem}\rm}
\newcommand{\mar}{\end{rem}}

\newcommand{\qd}{\end{proof}\vspace{0.5ex}}

\begin{document}

\title[$L^p$ 2-isometries]{A classification for 2-isometries of noncommutative
$L_p$-spaces}
\author{Marius Junge$^*$}
\address{Department of Mathematics\\
University of Illinois, Urbana, IL 61801, USA}
\email[Marius Junge]{junge@math.uiuc.edu}
\author{Zhong-Jin Ruan$^*$}
\email[Zhong-Jin Ruan]{ruan@math.uiuc.edu}
\author{David Sherman}
\email[David Sherman]{dasherma@math.uiuc.edu}
\thanks{${}^*$ The first and second authors were partially supported by the
National Science Foundation DMS-0088928 and DMS-0140067}

\begin{abstract}
In this paper we extend previous results of Banach, Lamperti and
Yeadon on isometries of $L_p$-spaces to the non-tracial case first introduced by Haagerup. Specifically, we use operator space techniques and an extrapolation argument to prove
that every 2-isometry $T : L_p(\M) \to L_p(\N)$ between arbitrary
noncommutative $L_p$-spaces can always be  written in the form
 \[ T(\p^{\frac 1p}) = w (\p\circ \pi^{-1} \circ E)^{\frac 1p}, \qquad \p \in \M_*^+.\]
Here  $\pi$ is a normal *-isomorphism
from $\M$  onto the von Neumann subalgebra $\pi(\M)$ of $\N$, $w$
is a partial isometry in $\N$,  and $E$ is a normal conditional
expectation from $\N$ onto $\pi(\M)$. As a consequence of this, any 2-isometry is automatically a complete isometry and has completely contractively complemented range.
\end{abstract}

\maketitle

\section{Introduction}

The investigation of isometries on $L_p$-spaces has a long
tradition in the theory of Banach spaces and has connections to
probability and ergodic theory. Banach \cite{Banach} considered
the discrete case and showed that for a surjective isometry on
$\ell_p$ one deduces the existence of  a permutation $\pi:
\mathbb{N} \to \mathbb{N}$ and a sequence of scalars
$\{\lambda_n\}$ with unit modulus so that $T(e_n) =\lambda_n
e_{\pi(n)}$.  Banach also stated, without proof, the form of
surjective isometries of $L_p([0,1], m)$. These are
\textit{weighted composition operators}; i.e.
\begin{equation} \label{E:banach}
T(f)(x) = h(x)f(\varphi(x)), \qquad x \in [0,1],
\end{equation}
where $\varphi$ is a measurable bijection of the unit interval
and $h$ is a measurable function with $|h|^p = \frac{d(m \circ
\varphi)}{dm}$. As shown by Lamperti \cite{Lamperti}, this
paradigm is still basically correct for non-surjective isometries
on general $L_p(X, \Sigma, \mu)$, but the ``composition" is
defined in terms of a mapping on the measurable sets (a so-called
\textit{regular set isomorphism}, see \cite{Lamperti} or
\cite{FJ}). For a sufficiently nice measure space, there is still
an underlying point mapping (see \cite{HvN}).  Notice that a set
mapping is nothing but a map on the projections in the associated
$L_\infty$-algebras, so Lamperti's result can be formulated
naturally in terms of von Neumann algebras!  Moreover, this
formulation implies that isometries must preserve
\textit{disjointness of support}, which has been a key step in
every stage of the classification of $L_p$-isometries, including
the present paper.

\bigskip

The most familiar noncommutative $L_p$-space is the Schatten-von
Neumann $p$-ideal $S_p$ \cite{vNSchatten}, the class of operators
in $\mathcal {K}(\mathfrak{H})$ (or in
$\mathcal{B}(\mathfrak{H})$) whose singular values are
$p$-summable. Segal \cite{Segal} extended this concept to general
semifinite von Neumann algebras. In this category, $L_p$-isometries have been
investigated by  Broise \cite {Broise}, Russo \cite {Russo},
Arazy \cite{Arazy}, Katavolos \cite {Kat1}, \cite {Kat2}, \cite
{Kat3}, and Tam \cite {Tam}. In 1981, Yeadon \cite{Yeadon} finally
gave a satisfactory  answer.

\bigskip

\begin{theo}
\label{P.Yeadon} \cite[Theorem 2]{Yeadon}
  Let $\RR$ and $\NN$ be two semifinite
von Neumann algebras with traces $\tau_\RR$ and $\tau_\NN$, respectively.
For $1\le p < \infty$ and $ p\neq 2$,  a
linear map
\[
T: L_p(\RR, \tau_\RR)
\to L_p(\NN, \tau_\NN) \]
is an isometry if and only if there exist a normal Jordan *-monomorphism
$J: \RR \to \NN $, a partial isometry $w \in \N$, and a positive self-adjoint
operator $B$ affiliated with $\N$ such that the spectral projections of $B$
commute with $J(\M)$,   which  verify the following conditions:
\begin{itemize}
\item [(1)] $w^* w = J(1) = s(B)$;
\item [(2)] $\tau_\M(x) = \tau_\N(B^p J(x))$ for all $x\in  \M^+$;
\item [(3)] $T(x) = w B J(x)$  for all $x \in  \M \cap L_p(\M, \tau_\M)$.
\end{itemize}
Moreover, $J$, $w$,  and $B$ are  uniquely determined
by the above conditions $(1) - (3)$.
\end{theo}


\bigskip

At approximately the same time, Haagerup \cite{Haagerup}
discovered a way to construct noncommutative $L_p$-spaces from
arbitrary von Neumann algebras. While many people worked on
general noncommutative $L_p$-spaces over the last twenty years, until recently little progress had been
made on the classification of their  isometries. Watanabe wrote a
series of papers (including \cite{W1} and \cite{W2}) supplying
many ideas and some partial results. In \cite{Sherman},
Watanabe's  techniques are developed to obtain classification
results  for $L_p$-isometries which are complementary to those
established in the present paper   (by entirely different
methods), and in \cite{Sherman2} a canonical form of surjective
$L^p$-isometries is derived, proving that $L_p(\mathcal{M})$ and
$L_p(\mathcal{N})$ are isometrically isomorphic if and only if
$\mathcal{M}$ and $\mathcal{N}$ are Jordan *-isomorphic.

\bigskip

In this paper, we use operator space techniques to give a complete classification of 2-isometries on general noncommutative $L_p$-spaces.  We begin by extending Yeadon's result to isometries
between noncommutative $L_p$-spaces where only the initial von
Neumann algebra is assumed  semifinite (using a recent result of
Raynaud/Xu \cite{RX}). After some background information (section
2) this task is performed in section 3. If we drop the assumption
that the initial algebra is semifinite, then $L_p(\M) \cap \M = \{0\}$, and a new kind of problem occurs. There is no longer a canonical embedding of the von
Neumann algebra (or finite elements therein) into the
noncommutative $L_p$-space, and therefore we may no longer use the
lattice of projections canonically. This has been the main
obstruction for applying Yeadon's groundbreaking technique
\cite{Yeadon}.  We found that operator space methods can be employed to
overcome this difficulty. The operator space structure of
noncommutative $L_p$-spaces was introduced by Pisier
\cite{PiLp} and later extended to the non-semifinite case (see
\cite{JRX} and \cite{Pibook}, the structure in \cite{Fidaleo} is
not entirely compatible here). In section 4 we are essentially
interested in taking columns and rows of elements in $L_p$ in order to classify \textit{complete} isometries.  Previous work for the case $p=1$ was done by \cite{Ozawa} and \cite{NO}. While isometries are connected with Jordan structure, a 2-isometry must already preserve multiplicative structure. This
is one of the key observations for our result, the other being an extrapolation argument which shows that if we have an $L_p$-isometry for one $p \neq 2$,
then we have an associated isometry for any other $p$. The case $p=4$ connects
with the theory of selfpolar forms (from \cite{Haa-St}) and thus can
be used to construct a conditional expectation, leading to a proof of the main theorems.  The final section of the paper consists of remarks.

To be more specific let us fix some notation. If $\N$ is a von Neumann
algebra, then $L_p(\N)$ is the linear span of elements
$\varphi^{\frac1p}$, where $\varphi$ ranges over the positive states
in $\N_*$.

\bigskip

\begin{theo}\label{main} Let $1\le p\neq 2<\infty$ and $T:L_p(\M)\to L_p(\N)$
be a linear map. The following are equivalent:
 \begin{enumerate}
 \item[i)] $id\ten T:L_p(M_2\ten \M)\to L_p(M_2\ten \N)$ is an
 isometry,
 \item[ii)] $id\ten T:L_p(M_m\ten \M)\to L_p(M_m\ten \N)$ is an
 isometry for all $m\in \nz$,
 \item[iii)] There exists an injective normal *-homomorphism $\pi:\M \to \N$,
 a partial isometry $w \in \N$, and a conditional expectation $E$ from $\N$ onto $\pi(\M)$ such that
  \begin{equation}\label{comp}  T(\varphi^{\frac1p})\lel w (\varphi\circ \pi^{-1}\circ E)^{\frac1p}, \qquad \varphi \in \N_*^+ \pl
  .
  \end{equation}
 \end{enumerate}
\end{theo}

In \eqref{comp} we have to use the fact that a map defined on
positive vectors extends uniquely to the whole of $L_p(\N)$.
Moreover, we recall that $L_p(\N)$ is an $\N-\N$ bimodule,
so that $w\varphi^{\frac1p}$ is a well-defined element satisfying
$\| w\varphi^{\frac1p}\|_p=\varphi(1)^{\frac1p}$ (see \cite{Terp1}
for more details). In the $\si$-finite case, $\phi^{\frac1p}\N$ is
dense in $L_p(\N)$ for every normal faithful state $\phi$. Then
we may alternatively describe $T$ by
\begin{equation}
T(\phi^{\frac1p}x)\lel w(\phi\circ
\pi^{-1}\circ E)^{\frac1p}\pi(x), \qquad x \in \M.
\end{equation}
Notice that when $\phi$ is not tracial, the inclusion mapping
$i_p:\N\to L_p(\N)$, $i_p(x)=\phi^{\frac1p}x$, no longer
preserves disjoint left supports. This is  the obstacle mentioned
above.

\bigskip

To end this introduction let us mention some recent sources which
will provide further background to interested readers.  The book
of Fleming and Jamison \cite{FJ} is an up-to-date survey on
isometries in Banach spaces (including classical $L_p$-spaces);
their companion volume will treat the case of noncommutative
$L_p$-spaces.  The handbook article of Pisier and Xu \cite{PX}
provides a general overview of noncommutative $L_p$-spaces,
focusing on Banach space properties and including a rich
bibliography. There is also recent groundbreaking work on the
\textit{non-isometric} isomorphism/embedding question by
Haagerup, Rosenthal,  and Sukochev \cite{HRS} which is entirely
disjoint to the isometric analysis of this paper.

\bigskip


\section{Some background}

In this section we provide some background on noncommutative $L_p$-spaces and
their operator space structure.
The readers  are referred to Nelson \cite{Nelson}, Kosaki \cite{Ko2},  Terp
\cite{Terp1} and  \cite {Terp2} for details on noncommutative $L_p$-spaces and
to Pisier \cite {PiLp}, \cite {Pibook}, Fidaleo \cite {Fidaleo}, and Junge-Ruan-Xu \cite{JRX} for the canonical
operator space structure on these spaces.
Pisier and Xu's recent survey paper  \cite{PX} provides a nice overview of this
subject.
For general background on the modular theory of von Neumann algebras,
we recommend the
Takesaki \textit{oeuvre} \cite{T1}, \cite{T2} and \cite{T3}, and
Str\v{a}til\v{a} \cite {Stratila}.

We first recall the $L_p$-space ($1 \le p < \infty$) associated with a
semifinite algebra $\mathcal{M}$ equipped with a given  normal faithful
semifinite trace $\tau$ (simply called a ``trace" from here on).
Consider the set
$$\{T \in \mathcal{M} \mid \|T\|_p \triangleq \tau(|T|^p)^{1/p} < \infty \}.$$
It can be shown that $\| \cdot \|_p$ defines a norm on this set.
Then the norm completion, which is denoted $L_p(\mathcal{M}, \tau)$,
is the noncommutative $L_p$-space obtained from $(\M, \tau)$.
It turns  out that one can identify elements of  $L_p(\mathcal{M}, \tau)$ with
certain $\tau$-measurable operators affiliated with $\M$
(see \cite{Nelson}).  Clearly $\tau$ is playing the role of integration here.

A von Neumann algebra lacking a faithful trace is not amenable to the previous
definition.
There are several alternative constructions which work in full generality.
Let us first recall the construction initiated by Haagerup \cite{Haagerup} and
carried out in detail by Terp \cite{Terp1}.
Choose a normal faithful semifinite weight $\phi$ on
$\mathcal{M}$. We consider the one-parameter modular automorphism $\sigma^\phi_t$
(associated with $\phi$) on $\M$ and   obtain a semifinite von
Neumann algebra
$\widetilde{\mathcal{M}} \triangleq \mathcal{M}
\rtimes_{\sigma^\phi} \mathbb{R}$
which  has an induced  trace  $\tau$ and a trace-scaling dual action
$\theta$ such that $\tau \circ \theta_s = e^{-s}\tau$ for all
$s \in \mathbb{R}$.

The original von Neumann algebra $\M$ can be identified with  a
$\theta$-invariant von
Neumann  subalgebra $L_\infty(\M)$ of $\widetilde{\mathcal{M}}$.
For $1 \le p < \infty$, the noncommutative $L_p$-space
$L_p(\M, \phi)$ is defined to be the space of all (unbounded) $\tau$-measurable
operators affiliated  with $\widetilde{\mathcal{M}}$ such that
$\theta_s(T) =  e^{-\frac sp}T$ for all $s \in \mathbb{R}$.
It is known from  Terp \cite[Chapter II]{Terp1} that there is a
one-to-one correspondence between bounded (positive) linear functionals
$\psi \in \M_*$ and
$\tau$-measurable (positive self-adjoint) operators $h_\psi
\in L_1(\RR, \phi)$ under the connection given by
\[
\widehat{\psi}(\tilde{x}) = \tau(h_\psi \tilde{x}), \qquad \tilde{x} \in \widetilde{\M},
\]
where $\widehat{\psi}$ is the so-called \textit{dual weight} for $\psi$.  This correspondence actually extends to all of $\M_*$ and $L_1(\RR, \phi)$, and we may define the ``tracial" linear functional $tr = tr_\M: L_1(\RR, \phi) \to \mathbb{C}$ by
\begin{equation}\label {F.conn}
tr(h_\psi) = \psi(1), \quad \text{satisfying} \quad tr (h_\psi x) = tr (x h_\psi) = \psi(x), \quad \psi \in \M_*^+, \: x \in \M.
\end{equation}

Given any $h \in L_p(\M)$, we have the polar decomposition $h = w |h|$,
where $|h|$ is a positive operator in $L_p(\M)^+$
and $w$ is a partial isometry contained in $\M$ such that the projection $s_l(h) = w w^*$ is the \emph{left support} of $h$ and
the projection $s_r(h) = w^* w$ is the \emph{right support} of $h$.  (We simply use $s$ for the support of positive vectors, operators, and maps.)
We can define  a Banach space norm on $L_p(\M, \phi)$ by
\begin{equation}
\label {F.lpnorm}
\|h\|_p = tr(|h|^p)^{\frac 1p} = \psi(1)^{\frac 1p}
\end{equation}
if $\psi\in \M_*$ corresponds to $|h|^p \in L_1(\M, \phi)^+$.
With this norm, it is easy to see that $L_1(\M, \phi)$ is
isometrically and orderly isomorphic to $\M_*$.
We note that  up to isometry
the noncommutative $L_p$-space constructed above is actually
independent of the choice of normal faithful semifinite weight on
$\M$. Therefore, we will simply write $L_p(\M)$
if there is no confusion.

We note that for any positive operator $h \in L_p(\M)^+$, $h^p$ is a positive
operator in $L_1(\M)^+$ and thus  we can write $h^p = h_\psi$ for a
corresponding  positive linear functional $\psi \in \M_*^+$.
Therefore, we may identify $h$ with $\psi^{\frac 1p}$,
i.e. we can simply write
\[
h = \psi^{\frac 1p}.
\]
This notation,  discussed specifically in \cite{Yamagami},
\cite[Section V.B.$\alpha$]{Connesbook} and  \cite{Sherman3}, provides the
relations
$$\psi^{it} \varphi^{-it}= (D\psi:D\varphi)_t, \qquad
\varphi^{it}x\varphi^{-it} = \sigma^\varphi_t(x), \qquad \varphi, \psi \in
\M_*^+, \; \varphi \text{ faithful}.$$
It will be used in section 4, where it suggests that the main results really deal with ``noncommutative
weighted composition operators".


For $1\le p < \infty$, we let $L_p(\M)'$ denote the dual space of $L_p(\M)$.
Then  we can obtain the isometric isomorphism
$L_p(\M)' = L_{p'}(\M)$ under  the \emph{trace duality}
\begin{equation}
\label {F.tracedual}
\langle x, y\rangle  = tr(xy) = tr(yx)
\end{equation}
for all $x \in L_p(\M)$ and $y \in L_{p'}(\M)$.  (Throughout $p'$ denotes the conjugate exponent of $p$; i.e. $\frac{1}{p} + \frac{1}{p'} = 1$.)

One of the advantages in using Haagerup's approach is that the natural
$\M-\M$ bimodule structure on $L_p(\M)$ is just operator composition.
It was shown in  \cite[Lemma 1.2]{JS} that for any $h \in L_p(\M)$, we have
\begin{equation} \label{E:cyclic}
\overline{\{x h \mid x \in \M\}} = L_p(\M)s_r(h) ~\text {and} ~
\overline{\{ h x \mid x \in \M\}} = s_l(h) L_p(\M).
\end{equation}
In particular, if $\mathcal{M}$ is a $\sigma$-finite von Neumann algebra and
$\phi$ is a normal faithful positive linear functional on $\M$,
then $h = \phi^{\frac 1p}$ is a \emph{cyclic vector} in $L_p(\M)$, i.e.
\[
\{x \phi^{\frac1p} \mid x \in \mathcal{M}\} ~ \text{and} ~
\{\phi^{\frac1p} x \mid x \in \mathcal{M}\}
\]
are norm dense in $L_p(\M)$.
In this case, we can obtain the following result
(see Junge and Sherman \cite[Lemma 1.3]{JS}).

\begin{lemma}
\label {P.conv}
Let $\phi \in \M_*^+$ be faithful.  A bounded net $\{x_\alpha\} \in \M$ converges strongly to $x$,
$x_\alpha  \overset{s}{\to} x$,  if and only if
$x_\alpha  \phi^{\frac1p} \to  x \phi^{\frac1p} $
(or $\phi^{\frac1p} x_\alpha   \to  \phi^{\frac1p} x$ )
in $L_p(\M)$.
\end{lemma}

Next let us recall Kosaki's complex interpolation construction of
$L^p(\mathcal{M})$  (see \cite{Ko2}).
We assume that $\M$ is a $\sigma$-finite von Neumann algebra and
$\phi$ is a  normal faithful state on $\M$.
Then we may identify $\M$ with a subspace of $\M_*$ by the
\emph{right embedding}
\[
x\in \M \mapsto \phi   x \in \M_*,
\]
where $\langle \phi  x, y\rangle = \phi(x y)$.
We use the right embedding (instead of the left embedding) in this paper
for the convenience of our  notation  in  representation theorems.
Then the complex interpolation $[\M, \M_*]_{\frac 1p}$ is a Banach space
which can be isometrically identified with
$\phi^{\frac{1}{p'}} L_p(\mathcal{M}) $ by
\[
\noo \phi^{\frac{1}{p'}} h \rrm_{[\M, \M_*]_{\frac 1p}} = \noo h\rrm_{L_p(\M)}
\]
for all $h \in L_p(\M)$ (see \cite[Theorem  9.1]{Ko2}).
In particular, any  $\phi^{\frac 1p} x \in L_p(\M)$ corresponds to
an element $\phi  x$ in $[\M, \M_*]_{\frac 1p}$ since
\begin{equation} \label{F.exchange}
\phi  x = \phi^{\frac{1}{p'}} (\phi^{\frac{1}{p}}x).
\end{equation}
In this case, we have
\begin{equation}
\label {F.norm}
\noo \phi x \rrm_{[\M, \M_*]_{\frac 1p}} = \noo \phi^{\frac 1p} x
\rrm_{L_p(\M)}.
\end{equation}
We will simply write $L_p(\M) = [\mathcal{M}, \mathcal{M}_*]_{\frac1p}$ when  there
is no confusion.

Using the complex interpolation, Pisier \cite {PiLp} contructed a
canonical operator
space matrix norm
\begin{equation}
M_n(L_p(\M)) = [M_n(\M), M_n(\M_*^{op})]_{\frac 1p}
\end{equation}
on $L_p(\M)$ (also see \cite  {Pibook} and \cite {JRX}).
For each $n\in {\Bbb N}$, we may also consider the noncommutative
$S^n_p$-integral
$S^n_p[L_p(\M)]$ of  $L_p(\M)$ and it turns out that we have the
isometric isomorphism
\begin{equation}
\label {Sp.norm} S_p^n[L_p(\RR)] =L_p(M_n\bar \otimes \RR).
\end{equation}
Moreover we can recover the canonical matrix norm on $L_p(\M)$ by
\begin{equation}
\label {F.Pisiernorm} \|x\|_{M_n(L_p(\M))}
= \mbox{sup}\{\|\alpha x \beta\|_{ S^n_p[L_p(\M)]}:
\|\alpha\|_{S^n_{2p}}, \|\beta\|_{S^n_{2p}} \le 1\}.
\end{equation}
Consequently, a linear map $T: L_p(\M) \to L_p(\N)$ is a complete contraction
(respectively, a complete isometry) if and only if  for
every $n \in {\Bbb N}$,
\[
\text{id}_{S^n_p} \otimes T : S^n_p[L_p(\M)] \to S^n_p[L_p(\N)]
\]
is a contraction (respectively, an isometry).


Finally let us recall the following equality condition for the noncommutative
Clarkson inequality, which  was shown for semifinite
von Neumann  algebras by Yeadon  \cite{Yeadon}, for general von Neumann
algebras with $2<p<\infty$ by Kosaki
\cite{Ko1}, and just recently for all $p \ne 2$ by Raynaud and Xu \cite{RX}.

\begin{theorem} \label{T:clarkson}
For $h, k \in L_p(\mathcal{M})$, $0 < p \neq 2 < \infty$,
\begin{equation}
\label {clak}
\|h + k\|_p^p + \|h - k\|_p^p = 2(\| h \|_p^p + \| k \|_p^p) \iff h
k^* = h^* k = 0.
\end{equation}
\end{theorem}

Let us agree to say that $L_p$ vectors $h$ and $k$ are \emph{orthogonal}
when they  satisfy the conditions in (\ref  {clak}).
We mentioned earlier that isometries on classical
$L_p$-spaces preserve disjointness of support.  Theorem \ref{T:clarkson} tells
us that isometries on general noncommutative $L_p$-spaces preserve
orthogonality.   This plays a key role in several of our proofs.

\section{A generalization of Yeadon's theorem} \label{S:gyt}

In 1981, Yeadon obtained a very satisfactory and complete
description for isometries between noncommutative $L_p$-spaces associated with
  semifinite von Neumann algebras (Theorem \ref{P.Yeadon} above).  We found that an analog of Yeadon's result still holds if the initial
algebra $\M$ is semifinite and the range algebra $\N$ is arbitrary.
The proof of this generalized result
is almost  the same as that given in Yeadon \cite {Yeadon}.
The major difference is that if $\N$ is a general von Neumann algebra
with a normal  faithful semifinite weight  $\phi$,
the   positive self-adjoint operator $B$ is affiliated with the
semifinite von Neumann algebra crossed product
$\widetilde \N = \N \rtimes_{\sigma^\phi} \mathbb{R}$ (instead of $\N$),
and condition (2) in Theorem \ref{P.Yeadon} should  be revised as
\[
\tau_\M(x) = tr_\N(B^p J(x))
\]
for all $x\in  \M^+$,
where $tr_\N$ is the Haagerup trace introduced in \eqref{F.conn}.
We will outline the proof of this generalized result in
the following theorem since it will provide us necessary notations
and motivation
for the rest  of  the paper.
We note that a proof by different methods can be found  in  Sherman
\cite{Sherman}.
Yeadon could not have proven such a result because general
$L_p$-spaces were only being invented as he was writing his paper, but he
\textit{did} prove it for preduals (the $p=1$ case) in \cite[$\S 4$]{Yeadon}.

In the rest of this section,  let us assume that  $\M$ is a semifinite
von  Neumann algebra with a (normal faithful semifinite) trace $\tau_\M$
and $\N$ is an arbitrary von Neumann algebra with a normal faithful semifinite
weight $\phi$, which induces the normal faithful semifnite trace
$\tau_{\widetilde \N}$
on the crossed product $\widetilde \N = \N \rtimes_{\sigma^\phi} \mathbb{R}$, and the Haagerup trace $tr_\N$ on $L_1(\N)$.

\begin{theorem}
\label{P.Jisometry}
For  $1\le p < \infty$ and $p \neq 2$, a  linear map
$$T: L_p(\M, \tau_\M) \to L_p(\N)$$ is an isometry if and only if there  exist a
normal  Jordan *-monomorphism  $J: \M \to \N$, a  partial isometry $w \in \N$,
and a  positive self-adjoint operator $B$ affiliated with $\widetilde \N$ such that $\theta_s(B) = e^{-\frac{s}{p}}B$ for all $s \in \mathbb{R}$ and the spectral projections of $B$ commute with $J(\M) \subset \N \subset \widetilde \N$,
which  verify the following
conditions:
\begin{enumerate}
\item [(1)] $w^* w = J(1) = s(B)$;
\item [(2)] $\tau_\M(x) = tr_\N(B^p J(x))$
for all $x\in  \M \cap L_p(\M)$;
\item [(3)] $T(x) = w  B  J(x)$ for all $x \in \M \cap L^p(\M, \tau_\M)$.
\end{enumerate}
Moreover, $J$, $w$,  and $B$ are  uniquely determined by
the above conditions $(1) - (3)$.
\end{theorem}
\begin{proof}
First note that the stated conditions do define an isometry; (2) implies
\[
\|x\|^p_p = \tau_\M( |x|^p) = tr_\N(B^p J(|x|^p)) =
tr_\N(|B J(x)|^p)  = \|T(x)\|^p_p
\]
for any $x \in \M\bigcap L_p(\M)$.  In the rest of the proof we derive these conditions for an arbitrary isometry $T$.

For each $\tau_\M$-finite projection $e\in \M$, we let $T(e) = w_e
B_e$ be the
polar decomposition  of $T(e)\in L_p(\N)$.   Then $w_e$ is a
partial isometry in $\N$ and $B_e = |T(e)|$ is a positive element of
$L_p(\N)$ such that
\begin{equation}
\label {F.commute}
B_e w_e^* w_e = B_e = w_e^* w_e  B_e.
\end{equation}
If we define $J(e) = w_e^* w_e = s_r(T(e)) = s(B_e)$ to be the corresponding
projection in $\N$, then $B_e$ commutes with $J(e)$.

Let $e$ and $f $ be two mutually orthogonal $\tau_\M$-finite
projections in $\M$.
Since $T$ is an isometry,  we have
\[
\|T(e) \pm T(f)\|_p^p = \|T(e \pm f)\|_p^p = \|e \pm f\|_p^p= \|e\|_p^p + \|f\|_p^p = \|T(e)\|_p^p +
\|T(f)\|_p^p.
\]
Then the Clarkson inequality is an equality, and we can
conclude by Theorem
\ref{T:clarkson} that
\[
T(e)^* T(f) = T(e) T(f)^* = 0.
\]
This implies that
\[
B_e B_f = 0, \qquad J(e)J(f) = 0, \qquad \text{and} \qquad w_e^* w_f = w_e w_f^* = 0.
\]
The linearity of $T$ gives $w_{e+f}B_{e+f} = w_e B_e + w_f B_f$, from which
$$B_{e+f}^2 = (w_{e+f} B_{e+f})^* (w_{e+f} B_{e+f}) = (w_e B_e + w_f B_f)^* (w_e B_e + w_f B_f) = B_e^2 + B_f^2 = (B_e + B_f)^2$$
and so
$$ w_e B_e + w_f B_f = w_{e+f}B_{e+f} = w_{e+f} B_e + w_{e+f} B_f.$$
This implies
\begin{equation}
\label {F.orth}
B_{e+f} = B_e + B_f, \qquad J(e+f) = J(e) + J(f), \qquad \text{and} \qquad  w_{e+f} = w_e + w_f.
\end{equation}

If $x = \sum \lambda_i e_i \in \M$ is a self-adjoint simple
operator with mutually orthogonal $\tau$-finite projections $\{e_i\}$  in $\M$,
we define $J(x) = \sum\lambda_i J(e_i)$.  It is easy to verify that
$$J(x^2) = J(x)^2 ~\text{and} ~  \|J(x)\|_\infty = \|x\|_\infty.
  $$
Moreover, we have $J(\lambda x) = \lambda x$  for all real $\lambda$
and $J(x + y) = J(x) + J(y)$ if $x$ and $y$ are  commuting
self-adjoint simple operators in $\M$.  In general, for any self-adjoint
operator $x \in \M$, there exists a sequence of simple functions $f_n$ on the
spectrum of $x$ with $f_n(0)=0$ which converges uniformly to the identity
function $f(\lambda) = \lambda$.  We define $J(x)$ to be the
$\|\cdot\|_\infty$-limit of $J(f_n(x))$ in $\N$.

Now fix a $\tau_\M$-finite projection $e \in \mathcal{M}$.  If $f$ is a projection in $\M$ such that $f \le e$, then $T(f) = T(e) J(f)$ and thus
\begin{equation}
\label  {F.repT}
  T(x) = T(e) J(x) = w_e B_e J(x)
\end{equation}
for  all $x \in e\M_{sa}e$.
It follows that for any $x, y \in e\M_{sa}e$,
\[
T(e) (J(x+y) - J(x) -J(y)) = T(x+y) - T(x) -T(y) = 0.
\]
This implies
\[
J(x+y) - J(x)  -J(y) = 0
\]
and thus $J$ is a real linear map from $e\M_{sa}e$ into $\N_{sa}$.
Next we can  extend $J$ to a *-linear map on all of $e\M e$ by
\[
J(x + iy) = J(x) + iJ(y), \qquad \forall x, y \in e\M_{sa}e.
\]
Still on $e\M e$, we have
\begin{eqnarray*}
J[(x + iy)^2] &=&  J[x^2 -y^2 + i( (x+y)^2 -x^2 - y^2)]  \\
&=&  J(x)^2 - J(y)^2 +  i((J(x)+ J(y))^2 -J(x)^2 - J(y)^2) = [J(x + iy)]^2.
\end{eqnarray*}
This shows that $J$ is  a Jordan  *-monomorphism on $e\M e$.
The normality of $J$ follows from \eqref {F.repT} and Lemma \ref  {P.conv}.
It is clear from \eqref {F.commute} and the above construction that all spectral
projections of  $B_e$ commute with  $J(e\M e)$.

If $\tau_\M$ is a finite
trace, we may take $B=B_1$ and the theorem is proved.
If $\tau_\M$ is not finite, some gluing must be done.
In this case, we may assume that  $\{e_\al\}$ is the (increasing) net of all
$\tau_\M$-finite projections in $\M$.
Then $e_\al \to 1$ in the strong operator topology.
By \eqref {F.orth} and \eqref  {F.repT},  we may  find a partial
isometry $w \in \N$ (as a strong limit of $\{w_{e_\al}\}$) and
a positive self-adjoint operator $B$ (as the supremum of $\{B_{e_\al}\}$)
affiliated with  $\widetilde \N$, and extend $J$ to a normal Jordan
*-monomorphism
from  $\M$ into $\N$ such that  conditions $(1)-(3)$ are satisfied.
In this case,  the spectral projections of $B$ commute with
$J(\M)$, and $B$ satisfies $\theta_s(B) = e^{-\frac sp} B$ for all $s
\in \Bbb R$.
We also have
\begin{equation}
\label {F.ind}
w_{e_\al} = wJ(e_\al),  ~ \, ~ B_{e_\al} =  B J(e_\al),
\end{equation}
and
\begin{equation}
\label {F.ind2}
tr_\N(B^p y) = \sup \{tr_\N(B_{e_\al}^p y)\}
\end{equation}
for every $y\in \N^+$.
Since each $B^p_{e_\al}\in L_1(\N)$ corresponds to a normal positive
linear functional
$\p_\al = tr_\N(B_{e_\al}^p \cdot)$ on $\N$,
$B^p$ corresponds to  a normal semifinite weight
$\p = tr_\N(B^p\cdot)$  on   $\N$, which has support $J(1)$.
Therefore, $\p$ is faithful when restricted to $J(1) \N J(1)$.

The uniqueness of $J,w,$ and $B$ is an easy consequence which we leave to the reader.
\qd

If we assume $J(1) = 1$ then $\p$ is a normal faithful semifinite weight on
$\N$. In this case, we may write $L_p(\N) = L_p(\N, \p)$ and it follows from
the conditions (2) and (3) in Theorem \ref {P.isometry} that the Jordan
*-monomorphism $J: \M \to \N$ induces an  orderly  isometric injection
$J_p^\p$ from  $L_p(\M, \tau_\M)$ into $L_p(\N, \p)$, which is given by
\begin{equation}
\label {F.inc}
J_p^\p(x) = B J(x)
\end{equation}
for all $x \in \M\bigcap L_p(\M, \tau)$.
If $J(1) \neq 1$, then we may orderly and isometrically identify
$L_p(\M, \tau_\M)$
with a subspace of
$L_p(J(1)\N J(1), \p) \simeq J(1)L_p(\N)J(1)$.

Using a standard matricial technique, we obtain the following
equivalence result.
The $p=1$ case is mentioned (without proof) in \cite{NO}.

\begin{prop}
\label{P.isometry}
Let $\M$ be a semifinite von Neumann algebra and let
$\N$  be an arbitrary von Neumann algebra.
For an isometry $T = w B J: L_p(\M, \tau_\M) \to L_p(\N)$
with $1\le p < \infty$ and $p \neq 2$, the following statements
are equivalent:
\begin{itemize}
\item [(i)] $T$ is a complete isometry,

\item [(ii)] $T$ is a $2$-isometry,

\item [(iii)] the Jordan map $J: \M \to \N$  is multiplicative.
\end{itemize}
\end{prop}
\begin{proof}
It is obvious that (i) $\Rightarrow$ (ii).

If $J$ is multiplicative from $\M$ into $\N$ then
$J_n= \text{id}_n \otimes J$  is a *-isomorphism from $M_n(\M)$ into
$M_n(\N)$ for every  $n \in \mathbb{N}$.  In this case, we can write
\[
id_{S^n_p} \otimes T : S^n_p [L_p(\M)] = L_p(M_n \bar \otimes
\M,  \text{tr}_n \otimes \tau_\M ) \to S^n_p[L_p(\N)] =L_p(M_n \bar \otimes
\N)
\]
as
\[
id_{S^n_p} \otimes T = w_n B_n J_n
\]
Here $w_n  = I_n \otimes w$ is a partial isometry
in $M_n(\N)$ such that $w_n^*w_n = J_n(1_n)$.  If we let $B = \sup \{B_{e_\alpha}\}$ as in the proof of Theorem \ref{P.Jisometry}, then
$B_n  = \sup_\alpha (B_{e_\alpha} \oplus_p B_{e_\alpha} \oplus_p \dots \oplus_pB_{e_\alpha})$ is a positive  selfadjoint
operator affiliated with $M_n(\widetilde \N)$ whose spectral projections commute with $J_n(M_n(\M))$.  We have that
\[
(\text{tr}_n \otimes \tau_\M) (x) = (\text{tr}_n \otimes
tr_\N)(B_n^p J_n(x))
\]
for all $x \in M_n(\M)^+$.
Then each $id_{S^n_p} \otimes T: S^n_p[L_p(\M)]  \to S^n_p[L_p(\N)] $
is an isometry by Theorem \ref {P.Jisometry}.
Therefore, $T$ is a complete isometry and we proved  (iii) $\Rightarrow$ (i).

It remains to prove (ii) implies (iii).
First $T$ is an isometry and thus has the representation $T = w B J$ by
Theorem \ref {P.Jisometry}.
Since  $T$ is also a $2$-isometry, we  deduce an isometry
\[
\widetilde T = {id}_{S^2_p} \otimes T: S_p^2[L_p(\M, \tau_\M)] \to
S_p^2[L_p(\N)] .
\]
Applying Theorem \ref{P.Jisometry} to $\widetilde T $, we
obtain a normal Jordan *-monomorphism
\[
\tilde J : M_2\bar \otimes \M \to M_2 \bar \otimes \N,
\]
a partial isometry $\widetilde w \in M_2 \bar \otimes \N$, and a
positive selfadjoint operator $\tilde B$
satisfying the conditions in Theorem \ref {P.Jisometry}.

If $e_i~ (i=1,2)$ are $\tau_\M$-finite projections in $\M$, then $\tilde e = \left[
\begin{matrix}
e_1& 0 \\
0 & e_2
\end{matrix}
\right]$ is a $(\text{tr}_2 \otimes \tau_\M)$-finite projection in $M_2 \bar \otimes \M$.  Let
$\widetilde T (\tilde e) =  w_{\tilde e}  B_{\tilde e}$ be the polar
decomposition of  $T_2(\tilde e)$  and let $T(e_i) = w_{e_i}
B_{e_i}$ be the
polar  decomposition of
$T(e_i)~ (i=1,2)$.  Since
\[
\widetilde T (\tilde e)  = \left[
\begin{matrix}
T(e_1)& 0 \\
0 & T(e_2)
\end{matrix}
\right] =
\left[
\begin{matrix}
w_{e_1}& 0 \\
0 & w_{e_2}
\end{matrix}
\right]
\left[
\begin{matrix}
B_{e_1} & 0 \\
0 & B_{e_2}
\end{matrix}
\right]
\]
we must have, by the uniqueness of the polar decomposition,
\[
w_{\tilde e} =
\left[
\begin{matrix}
w_{e_1}& 0 \\
0 & w_{e_2}
\end{matrix}
\right]
~ \mbox{and} ~
B_{\tilde e}  = \left[
\begin{matrix}
B_{e_1} & 0 \\
0 & B_{e_2}
\end{matrix}
\right].
\]
According to the definition of $\widetilde J $ given in the proof of
Theorem \ref {P.Jisometry}, we have
\[
\widetilde J \left (\left[
\begin{matrix}
e_1& 0 \\
0 & e_2
\end{matrix}
\right]\right ) = \left[
\begin{matrix}
J(e_1)& 0 \\
0 & J(e_2)
\end{matrix}
\right],
\]
and thus we can conclude that
\[
\widetilde J \left (\left[
\begin{matrix}
x& 0 \\
0 & y
\end{matrix}
\right]\right ) = \left[
\begin{matrix}
J(x)& 0 \\
0 & J(y)
\end{matrix}
\right]
\]
for all $x, y \in \M$.  Since $\widetilde T$ is an $M_2$-bimodule morphism, we can write
\[
\widetilde T \left ( \left[
\begin{matrix}
0& x \\
y & 0
\end{matrix}
\right]\right )
= \widetilde T \left ( \left[
\begin{matrix}
   0& 1 \\
1 & 0
\end{matrix}
\right]
\left[
\begin{matrix}
x &0 \\
0& y
\end{matrix}
\right]\right )  =
\left[
\begin{matrix}
0& 1 \\
1 & 0
\end{matrix}
\right]
\left[
\begin{matrix}
T(x) &0 \\
0& T(y)
\end{matrix}
\right].
\]
Then for any $x, y \in \M$, we obtain
\[
\widetilde J\left (
\left[
\begin{matrix}
0& x \\
y & 0
\end{matrix}
\right] \right ) =
\left[
\begin{matrix}
0& J(x) \\
J(y) & 0
\end{matrix}
\right].
\]
 From this we can conclude that $J(xy) = J(x)J(y)$  since
\begin{eqnarray*}
\left[
\begin{matrix}
J(xy) & 0 \\
0 & J(yx)
\end{matrix}
\right] &=&
\widetilde J \left (\left[
\begin{matrix}
xy & 0 \\
0 & yx
\end{matrix}
\right] \right )
=  \widetilde J \left (\left[
\begin{matrix}
0 & x\\
y & 0
\end{matrix}
\right]^2  \right ) =
  \widetilde  J  \left (\left[
\begin{matrix}
0 & x\\
y & 0
\end{matrix}
\right]  \right )^2 \\
&=&
\left [  \begin{matrix}
0 & J(x) \\
J(y) & 0
\end{matrix}
\right]
\left [ \begin{matrix}
0 & J(x) \\
J(y) & 0
\end{matrix}
\right] =
\left [\begin{matrix}
J(x)J(y) &0\\
0 & J(y)J(x)
\end{matrix}
\right].
\end{eqnarray*}
Therefore $J$ is multiplicative.
Finally we note that we can conclude from the above calculations that
$\widetilde J = J_2$, $\widetilde w = w_2$ and
$\widetilde B = B_2$.
\qd

If $T = w B J: L_p(\M, \tau_\M) \to L_p(\N)$ is a $2$-isometry (or equivalently,
a complete isometry) with $J$ a normal *-monomorphism from $\M$ into
$\N$,  then
\begin{equation}
\label {F.smap1}
S = w^* T = B J
\end{equation}
is a completely positive and completely isometric injection from
$L_p(\M, \tau_\M)$ into  $L_p(\N)$.
Therefore, we can completely orderly and completely isometrically
identify $L_p(\M, \tau_\M)$ with the operator subspace $S(L_p(\M, \tau_\M))$ of $L_p(\N)$.

\begin{prop} \label{T:complemented}
Let $T: L_p(\M, \tau_\M) \to L_p(\N)$ be a $2$-isometry (or equivalently, a complete isometry).
Then $T(L_p(\M, \tau_\M))$ is completely contractively complemented
in $L_p(\N)$.

If, in addition, $T$ is positive, then $T(L_p(\M, \tau_\M))$ is completely
positively and
completely contractively complemented  in   $L_p(\N)$.
\end{prop}

\begin{proof}
Assume the notation of Theorem \ref{P.Jisometry}.  We have that $J$ is multiplicative by Proposition \ref{P.isometry}, so that $J(\M)$ is a von Neumann subalgebra of $\N$.  Without loss of generality, we may assume that $\p = tr_\N(B^p \cdot)$ is faithful on $\N$.
Otherwise we restrict our argument to $J(1) \N J(1)$.

Since $B^p$ is the Pedersen-Takesaki derivative \cite {PT} of $\widehat{\varphi}$ on $\widetilde{\N}$ with respect to $\tau_{\widetilde{\N}}$, we have that
$$\sigma_t^\p(y) = \sigma_t^{\widehat{\p}} (y) = B^{ipt}y B^{-ipt}, \qquad y \in \N \subset \widetilde{\N}.$$
In particular, we have
\[
\sigma^{\p}_t (J(x))= B^{ipt}J(x) B^{-ipt} = J(x)
\]
for all $x \in \M$, since $B$ commutes with $J(\M)$.
Then we can conclude from Takesaki's theorem \cite {Ta} that
there exists a unique normal conditional expectation $E$ from
$\N$ onto $J(\M)$ such that
\[
\p\circ E = \p.
\]
 From this, we may induce a completely positive and completely contractive projection $E_p$ from
$L_p(\N)= L_p(\N, \p)$ onto $w^*T(L_p(\M, \tau))$ (see \cite{HRS}, \cite {JX}, \cite{Sherman}, or the construction sketched after Lemma \ref{41}).  Then $h \mapsto w E_p(w^* h)$ is the required projection.  If $T$ is positive, $w$ and $w^*$ may be omitted, so this projection is $E_p$ itself.
\qd

Proposition \ref{T:complemented} is a noncommutative version of the classical fact that for any isometric
embedding
$$T: L_p(X, \Sigma_X, \mu_X) \to L_p(Y, \Sigma_Y, \mu_Y), $$
the image space $T(L_p(X, \Sigma_X, \mu_X))$ must be contractively
complemented in
$L_p(Y, \Sigma_Y, \mu_Y)$ (see Lacey \cite{Lacey}).  Stronger results are in Theorem \ref{T:sigmafinite} of this paper and Section 8 of \cite{Sherman}.

\section{2-isometries on general noncommutative $L_p$-spaces} \label{S:2isom}
\newcommand{\bfi}{\bar{\phi}}

The aim of this section is to study 2-isometries on
general noncommutative $L_p$-spaces.  We will use the alternative notation,
i.e. positive linear  functionals $\p^{\frac 1p}$, instead of positive
self-adjoint operators $h_\p^{\frac 1p}$, for elements in $L_p(\M)^+$.  Until the proof of Theorem \ref{main} given at the end of the section, we assume that $\M$ is an arbitrary $\sigma$-finite  von Neumann algebra with a fixed normal faithful state $\phi$.  We also use $\phi$ for its corresponding density operator $h_\phi$ in $L_1(\M)$, so that $\phi^{\frac 1p} \M$ is norm dense in $L_p(\M)$.

As in \S 3, we start by using support projections to construct an embedding from $\M$ into $\N$.  But since $\M \bigcap L_p(\M) = \{0\}$ ($p \ne \infty$) when $\M$ is not semifinite, we cannot directly employ the projection lattice of $\M$, and instead work with the vectors $\{\p^{\frac1p} x \mid x \in \M\}$.  Then $2\times2$ matrix equations produce the *-monomorphism $\pi$ (in the next proposition), but several steps are still required to show that $\pi(\M)$ is complemented and obtain the desired decomposition for $T$.


\begin{prop} \label{T:mod}
Let $1\le p<\infty$ and $p\neq 2$. If $T:L_p(\M)\to
L_p(\N)$ is a $2$-isometry, then
there exists a normal *-monomorphism $\pi:\M\to \N$ such that
\begin{equation}
\label {sta}
T(\phi^{\frac 1p} x)\lel T(\phi^{\frac 1p}) \pi(x) \pl, \qquad x \in \M.
\end{equation}
Moreover, $\pi$ does not depend on the choice of (faithful) $\phi \in \M_*^+.$
\end{prop}
\begin{proof}
We first define a map $\pi$ between projection lattices by
\begin{equation}
\label{E:proj}
\pi: e \in \P(\M)  \mapsto s_r(T(\phi^{\frac1p} e )) \in \P(\N).
\end{equation}
That is, with  the polar  decomposition
$T(\phi^{\frac1p }e) = w_e |T(\phi^{\frac1p }e)|$,
we set $\pi(e) = w_e^* w_e $.
If $e \perp f$, then
$$ \left[ \begin{matrix} \phi^{\frac1p} e & 0 \\ 0 & 0 \end{matrix}
\right] ~ \mbox{and} ~  \left[ \begin{matrix} 0 & 0 \\
\phi^{\frac1p} f & 0
\end{matrix} \right ]$$
are two orthogonal elements in $S^2_p[ L_p(\M)]$.
By Theorem \ref{T:clarkson}, their images
$$ \left [ \begin{matrix} T(\phi^{\frac1p} e) & 0 \\ 0 & 0 \end{matrix}
\right ] ~ \mbox{and} ~
\left [ \begin{matrix} 0 & 0 \\ T(\phi^{\frac1p} f) & 0
\end{matrix} \right ]$$
under $\text{id}_{S^2_p}  \otimes T$  are orthogonal in $S^2_p[ L_p(\N)]$.
This shows that
$$ s_r \left (\left [ \begin{matrix} T(\phi^{\frac1p} e) & 0 \\ 0 & 0
\end{matrix}  \right ]\right ) = \left [ \begin{matrix}
s_r(T(\phi^{\frac1p} e)) & 0 \\ 0 & 0  \end{matrix} \right ]
= \left [ \begin{matrix}
\pi(e) & 0 \\ 0 & 0  \end{matrix} \right ]$$
is orthogonal to
$$ s_r \left (\left [ \begin{matrix} 0 & 0 \\ T(\phi^{\frac1p} f) & 0
\end{matrix}  \right ]\right ) = \left [ \begin{matrix}
s_r(T(\phi^{\frac1p} f)) & 0 \\ 0 & 0  \end{matrix} \right ]
= \left [ \begin{matrix}
\pi(f) & 0 \\ 0 & 0  \end{matrix} \right ],
$$
and thus $\pi(e)$ is orthogonal to $\pi(f)$.
Since
$$T(\phi^{\frac1p})\pi(e) = [T(\phi^{\frac1p} e) + T(\phi^{\frac1p}(1-e))]
\pi(e)  = T(\phi^{\frac1p }e) + T( \phi^{\frac1p}) \pi(1-e) \pi(e) =
T(\phi^{\frac1p} e),$$
we have \eqref {sta} for projections in $\M$.
As a consequence, we obtain
\[
\pi(e+f) = \pi(e) + \pi(f)
\]
for orthogonal projections $e\perp f$ since
$$T(\phi^{\frac1p})\pi(e+f)  = T(\phi^{\frac1p}(e+f)) = T(\phi^{\frac1p} e) +
T(\phi^{\frac1p} f) = T(\phi^{\frac1p})(\pi(e) + \pi(f))$$
and $T(\phi^{\frac1p})$ is separating for the right action of $\pi(1)\N\pi(1)$.

Now we extend $\pi$ as in the proof of Theorem \ref{P.Jisometry}: first to
finite real linear combinations of orthogonal projections, then to all
self-adjoint elements in $\M$ by continuity, and finally to all of $\M$ by complex linearity.  Apparently $\pi$ satisfies
$$ T(\phi^{\frac1p}x) = T(\phi^{\frac1p} ) \pi(x), \qquad x \in \M.$$
This relation implies additivity: for $x,y \in \M$, we have
$$T(\phi^{\frac1p})\pi(x + y)  = T(\phi^{\frac1p}(x+y)) = T(\phi^{\frac1p }x)
+ T(\phi^{\frac1p}y) = T(\phi^{\frac1p})(\pi(x) + \pi(y)).$$

Now let $u$ be a unitary element of $\M$.  Replacing $\phi^{\frac1p}$ by $
\phi^{\frac1p}u$ in the definition of $\pi$, we obtain a new *-preserving map
$\pi_u$.  For a projection $e$ in $\M$,
$$ \left [ \begin{matrix} \phi^{\frac1p} e & 0 \\ 0 & 0 \end{matrix}
\right ] \text{ and } \left [\begin{matrix} 0 & 0 \\ \phi^{\frac1p} u(1-e)
& 0 \end{matrix} \right ]$$
are orthogonal, so $\pi(e) \perp \pi_u(1-e)$.  Since
$$\pi_u(1) = s_r(T(\phi^{\frac1p}u)) = s_r( \pi(u)) = \pi(1),$$
we must have $\pi(e) = \pi_u(e)$ for every projection $e$, whence $\pi =
\pi_u$.  Then for any $x \in \M$,
$$T(\phi^{\frac1p})\pi(ux)  = T( \phi^{\frac1p}ux) = T(\phi^{\frac1p} u)\pi(x)
  = T(\phi^{\frac1p})\pi(u) \pi(x).$$
Since the linear span of the unitaries is all of $\M$, it follows that $\pi$
is multiplicative.  Then
$$T(\phi^{\frac1p} xy) = T(\phi^{\frac1p}) \pi(xy) =
T(\phi^{\frac1p})\pi(x)\pi(y) = T(\phi^{\frac1p}x)\pi(y), \qquad x,y \in \M,$$
(densely) establishes \eqref{sta} and the independence of $\pi$ from the choice of $\phi$.

To prove the normality of $\pi$, we let  $\{x_\alpha\}$ be a bounded
net in $\M$
such that  $\notag x_\alpha \overset{s}{\to} x$ in the strong topology.
It follows from Lemma \ref {P.conv} that  $\phi^{\frac1p} x_\alpha \to
\phi^{\frac1p} x$ in $L_p(\M)$, and thus
\[
T(\phi^{\frac1p}) \pi(x_\a) = T(\phi^{\frac1p}x_\alpha )
\to T(\phi^{\frac1p}x) = T(\phi^{\frac1p})\pi(x)
\]
in $L_p(\N)$.
Again by Lemma \ref {P.conv}, this implies
\[
\pi(x_\alpha) \overset{s}{\to} \pi(x)
\]
in  $\pi(1) \N \pi(1)$ and thus in $\N$.
\end{proof}

Duality will be a very important tool in the following arguments.
Given a bounded linear map $T:L_p(\M)\to L_p(\N)$, we let
$T': L_{p'}(\N)\to L_{p'}(\M)$ denote the adjoint of $T$ and let
$T^{'*} = (T')^*$ denote the *-adjoint of
$T'$, i.e.,
\[
T^{'*}(k) = T'(k^*)^* \qquad k\in
L_p(\N).
\]
Then
\[
tr_\M(T^{'*}(k)^*h) \lel tr_\N(k^* T(h)) \pl, \qquad h\in
L_p(\M), \: k \in L_{p'}(\N),
\]
defines a sesquilinear form on $L_p(\M) \times L_{p'}(\N)$.
We also set the notation $\bar{\varphi} \triangleq |T(\varphi^{\frac1p})|^p$ for any state $\varphi \in \M_*^+$.  Thus $\bar{\varphi}^{\frac1p}$ is the absolute
value of $T(\varphi^{\frac1p} )$.

\begin{lemma} \label{comp1} Let $1\le p<\infty$  and
let $T:L_p(\M)\to L_p(\N)$ be a
$2$-isometry (or an isometry satisfying \eqref{sta}). Then
\[
\bar{\phi}\circ \pi \lel \phi \pl .
\]
\end{lemma}
\begin{proof} Let us first assume that  $p \ne 1$.
Since $\phi$ is a normal faithful state in $\M_*^+$, then
$\phi^{\frac 1p}$ is a unit
vector in $L_p(\M)$ since
\[
\|\phi^{\frac 1p}\|_p ^p = \phi(1) = 1.
\]
This implies that $T(\phi^{\frac1p})=w\bar{\phi}^{\frac1p}$ is a unit vector in
$L_p(\N)$ and   thus
  \[ 1\lel tr_\N((w\bar{\phi}^{\frac{1}{p'}})^*T(\phi^{\frac1p}))
   \lel tr_\M(T^{'*}(w\bar{\phi}^{\frac{1}{p'}})^*\phi^{\frac1p}) \pl .\]
Since $L_{p'}(\M)$ is uniformly convex,
$\phi^{\frac1p}$ admits exactly one norm attaining element. Thus we must have
\[ T^{'*}(w \bfi^{\frac{1}{p'}})\lel \phi^{\frac{1}{p'}} \pl .\]
On the other hand, we have
  \begin{eqnarray*}
   \phi(x)&=& tr_\M((\phi^{\frac{1}{p'}})^* \phi^{\frac1p}x) \lel
   tr_\M(T^{'*}(w \bfi^{\frac{1}{p'}})^* \phi^{\frac1p}x) \lel
tr_\N((w \bfi^{\frac{1}{p'}})^* T(\phi^{\frac1p}x)) \\
&\lel&  tr_\N(\bfi^{\frac{1}{p'}} w^* w \bfi^{\frac1p} \pi(x))
\lel \bfi(\pi(x)) \pl.
   \end{eqnarray*}
If $p=1$, the same argument applies, replacing $\bfi^{\frac{1}{p'}}$ and
$\phi^{\frac{1}{p'}}$ by 1.  (Since $\phi$ is faithful, it attains its norm at
1 only.)
  \qd

We can already say quite a lot in case $p=1$.

\begin{prop} \label{T:lone}
Let $T: \M_* \to \N_*$ be a $2$-isometry.  There are  a normal
*-monomorphism $\pi: \M \to \N$, a partial isometry $w \in  \N$,
and a normal conditional  expectation $E:\N \to \pi(\M)$ such that
$$T(\p) = w(\p \circ \pi^{-1} \circ E), \qquad \p \in \M_*.$$
\end{prop}
\begin{proof}
It was shown in \cite[Theorem 3.2]{Sherman} that
for any such $L_1$-isometry $T$, there
are a normal Jordan *-monomorphism $J: \M \to \N$, a partial isometry
$w \in \N$, and a normal positive projection $P:\N \to J(\M)$, faithful
on $J(1)\N  J(1)$, such that
$$T(\p) = w(\p \circ J^{-1} \circ P) , \qquad \p \in \M_*.$$
(We note that an equivalent result was proved earlier by Kirchberg
\cite{Kirchberg}.)  So the induced map $S = w^* T$ is
a  positive isometry and thus we can get $\bar \p = \p \circ J^{-1} \circ P$
for  all $\p \in \M_*^+$.  Precomposing with the $\pi$ from Proposition \ref{T:mod} and applying Lemma \ref{comp1},
$$\p \circ J^{-1} \circ P \circ \pi = \bar \p \circ \pi = \p, \qquad \forall \text{ faithful } \p \in \M_*^+.$$
Apparently $ J^{-1} \circ P \circ \pi$ is the identity map.  In particular,
we have  $P(\pi(e)) = J(e)$ for any projection $e$ in $\M$.
Using properties of $P$ and $J$  (see \cite[Lemma 5.4]{Sherman}),
we can conclude that
$$P(J(e)\pi(e)J(e)) = J(e)P(\pi(e))J(e) = [J(e)]^3 = J(e) = P(J(e)).$$
Since  $J(e)\pi(e)J(e) \le J(e)$ and $J(1) = s(P)$, we must have
$J(e)\pi(e)J(e) = J(e)$.
This implies $\pi(e) \ge J(e)$.  Since  $P(\pi(e) - J(e)) = 0$,  we must
have $\pi(e) = J(e)$.  Therefore $J=\pi$ is multiplicative and $E=P$ is a
normal conditional expectation from $\N$ onto $\pi(\M)$.
\qd

\begin{lemma}\label{intpol} Let $\N_1$ be a von Neumann subalgebra of $\N_2$
and  let $\bfi$ be a normal faithful state on $\N_2$ with $\phi=\bfi|_{\N_1}$.
For $2\le p\le \infty$,
  \[ \noo \bfi^{\frac1p} x\rrm_{L_p(\N_2)}\kl  \noo
  \phi^{\frac1p}x \rrm_{L_p(\N_1)} \]
for all $x\in \N_1$.
\end{lemma}
\begin{proof}
We first recall from Section 2 that we can identify  Haagerup's
$L_p$-space  $L_p(\N_i)$ with the  complex interpolation spaces
$[\N_i, (\N_i)_*]_{\frac 1p}$ ($i=1,2$) (see \eqref {F.exchange}
and \eqref{F.norm}).
Since $\phi=\bfi|_{\N_1}$ and $\bfi$ are  normal faithful states on $\N_1$ and
$\N_2$, respectively,
the induced inclusion
\[
\iota_2 : \phi^{\frac 12} x \in L_2(\N_1, \phi) \to \bfi^{\frac 12} x
\in L_2(\N_2,
\bfi)
\]
is an isometric inclusion.
The corresponding inclusion
\[
\iota_2: \phi x \in [\N_1, (\N_1)_*]_{\frac 12} \to \bfi x \in
[\N_2, (\N_2)_*]_{\frac 12}\]
between complex interpolation spaces  is also an isometric inclusion.
Then for any $2\le p < \infty$ the canonical inclusion
\[
\iota_p: \phi ^{\frac 1p} x \in L_p(\N_1) \to \bfi^{\frac 1p} x \in L_p(\N_2)
\]
is a contraction  since it can be identified with the complex
interpolation
\[
\iota_p = [\iota_\infty, \iota_2]_{\frac 2p},
\]
where we let $\iota_\infty : \N_1 \hookrightarrow \N_2$   denote the
canonical  inclusion  of $\N_1$ into $\N_2$.
\qd

Getting back to $T$, let us assume that
\[
T(\phi^{\frac 1p} x) = \bfi^{\frac 1p} \pi(x).
\]
Otherwise, we may  replace $T$ with
$\phi^{\frac 1p} x  \mapsto  w^* T(\phi^{\frac 1p} x)$
where $w$ is the partial isometry obtain from the
polar decomposition $T(\phi^{\frac1p}) = w \bfi^{\frac1p}$.
For any $1\le q < \infty$, we define the related maps
  \[
T_q(\phi^{\frac1q}x)\lel \bar{\phi}^{\frac1q} \pi(x)\pl .\]

\begin{cor}\label{dual1} Let $1<p\le 2$.
If $T = T_p: \phi^{\frac1p} x \mapsto \bfi^{\frac1p}\pi(x)$ is an isometry, then $T_{p'}$ is also an isometry.
\end{cor}
\begin{proof}
Let us show that  $T_{p'}^{'*}T_p$ is the identity on $L_p(\M)$.
Indeed, by definition and Lemma \ref{comp1} we find
  \begin{align*}
   tr_\M( T_{p'}^{'*}(\bfi^{\frac1p}\pi(y))^*\phi^{\frac{1}{p'}}x)
   &= tr_\M((\bfi^{\frac1p}\pi(y))^*T_{p'}(\phi^{\frac{1}{p'}}x))
   = tr_\N(\pi(y^*) \bfi^{\frac{1}{p}}\bfi^{\frac{1}{p'}}\pi(x)) \\
   &\lel \bfi(\pi(xy^*)) \lel \phi(x y^*)
   = tr_\M(( \phi^{\frac{1}{p}}y)^*\phi^{\frac{1}{p'}}x) \pl .
  \end{align*}
By duality we conclude that $T_{p'}^{'*}(\bfi^{\frac1p}\pi(y)) = \phi^{\frac{1}{p}}y$.  This implies that
\[
T_{p'}^{'*}(T_{p}(\phi^{\frac1p}y))) =
T_{p'}^{'*}(\bfi^{\frac1p}\pi(x))=\phi^{\frac{1}{p}}y.
\]
This shows that  $T_{p'}^{'*}T_p = {\rm id}_{L_p(\M)}$.
By duality, we deduce that
  \[
\text{id}_{L_{p'}(\M)}\lel
(T_{p'}^{'*}T_p)^{'*} \lel T_p^{'*}T_{p'} \pl .\]
By assumption $T_p$ and thus $T_p^{'*}$ is a contraction.  $T_{p'}$ is also
  contractive by Lemma \ref{intpol}, so it must be an  isometry.
\qd

Due to Corollary \ref{dual1}, we may now focus on the case
$2<p<\infty$. The main argument here is to show that if $T_p$ is
isometric for some value of $p$, then it is isometric for all
values. We will use Kosaki's interpolation theorem for the proof
of this result but in a more explicit form. Let us use the
notation $\mathfrak{S} =\{ z\in \cz \pl|\pl 0\le Re(z)\le 1\}$.

\begin{lemma}\label{exKos} Let $\phi$ and $\psi$ be normal faithful states
such that $\psi \le C \phi$ for some constant $C>0$. Let
$2< r< p< q < \infty$ and
$\frac{1}{p}=\frac{1-\theta}{q}+\frac{\theta}{r}$. Then there
exists an analytic function $h: \mathfrak{S} \to \M$ such that
   \begin{enumerate}
  \item $\noo \phi^{\frac{1}{q}} h(it) \rrm_q\le 1$,
   \item $\noo \phi^{\frac{1}{r}} h(1+it) \rrm_r\le 1$,
   \item $ \phi^{\frac1q} h(0) =\psi^{\frac1q}$, $ \phi^{\frac1p} h(\theta)
\lel \psi^{\frac1p}$,
   $ \phi^{\frac1r} h(1) =\psi^{\frac1r}$,
  \item the maps $t\mapsto h(1+it)\phi^{\frac1r}$, $t\mapsto
h(it)\phi^{\frac1q}$ are continuous.
  \end{enumerate}
\end{lemma}
\begin{proof}
Since $\psi \le C \phi$, we have by \cite[Theorem VIII.3.17]{T2} (taking
adjoints) that the Connes cocycle derivative $(D\phi: D\psi)_t = \psi^{it}
\phi^{-it}$ extends off the real line to a $\sigma$-weakly continuous function
on $\{z\,\mid \, 0 \le \text{Im}\, z \le 1/2 \}$ which is analytic in the
interior.  We define $\frac{1}{s}\lel \frac{1}{r}-\frac{1}{q}$ and set
$$h(z) = (D\phi:D\psi)_{i\left(\frac{1}{q} + \frac{z}{s} \right)} =
\phi^{-\frac{1}{q}-\frac{z}{s}}\psi^{\frac{1}{q}+\frac{z}{s}}.$$
Note that $h$ is analytic on $\mathfrak{S}$.

We have
  \begin{align*}
   \noo \phi^{\frac1q} h(it) \rrm &=  \noo
   (\phi^{-\frac{it}{s}}\psi^{\frac{it}{s}}) \psi^{\frac1q} \rrm_q
   \lel \noo \psi^{\frac1q} \rrm_q
    \lel 1. \pl \
  \end{align*}
Similarly,
  \[ \noo \phi^{\frac1r}  h(1+it) \rrm_r \lel \noo
  (\phi^{-\frac{it}{s}}\psi^{\frac{it}{s}})  \psi^{\frac1r} \rrm_q
\lel 1 \pl .\]
The equalities for $ \phi^{\frac1q} h(0) $, $ \phi^{\frac1r} h(1) $ and
$ \phi^{\frac1p} h(\theta) $ are obvious.  Since $(
\phi^{-\frac{it}{s}}\psi^{\frac{it}{s}})$ is a cocycle, it is strongly
continuous. Then we obtain the last  statement by
Lemma \ref {P.conv} and the above equations.
\qd

We are now able to prove the key extrapolation result.

\begin{prop}\label{extra} Let $\phi$ be a  normal faithful state on $\N$,
$\bar{\phi}\circ
\pi=\phi$, and  $T_p(\phi^{\frac1p}x)\lel
\bar{\phi}^{\frac1p}\pi(x)$. If $T_p$ is an isometry for some
$2<p<\infty$, then $T_p$ is an isometry for all $2<p<\infty$.
\end{prop}
\begin{proof} Let $2< r<p<q<\infty$ and assume that $T_p$ is an isometry.  We
want to show that $T_q$ and $T_r$ are isometries.  We define $\frac1s =
\frac{1}{r}-\frac{1}{q}$ and $\theta$ by
$\frac{1}{p}=\frac{1-\theta}{q}+\frac{\theta}{r}$.  Let us assume that $\psi$
is a normal faithful  state with $\psi\le C\phi$ for some constant $C <
\infty$.  (Such $\psi$ form a dense face of $\M_*^+$, see e.g. \cite{Ju}.)  We
know that $T_q$ is a contraction by Lemma \ref{intpol}; assume toward a
contradiction that
  \[ \noo T_q(\psi^{\frac1q}) \rrm_q < 1 \pl .\]
Since $T_q$ is continuous, we can find a $\delta>0$ such that
\begin{equation} \label{E:notisom} \noo T_q( \phi^{-\frac{it}{s}}
\psi^{\frac{it}{s}} \psi^{\frac1q})\rrm_q\le
  (1-\delta)
\end{equation}
for all $|t|\le \delta$. Let $\mu_\theta$ be the probability
measure on the boundary of the strip $\mathfrak{S}$ such that
  \[ \intt_{\partial \mathfrak{S} } f(z) d\mu_\theta(z)\lel f(\theta) \]
for every harmonic function on $ \mathfrak{S} $.  This is just a relocation of
the Poisson kernel, so Lebesgue measure is absolute continuous with respect to
$\mu_{\theta}$ (an explicit formula is in \cite[Section 4.3, p.93]{BL}).
Therefore $\mu_{\theta}([-i\delta,i\delta])>0$.  So we may also find $\eps>0$
such that $\frac{1-\eps}{1+\eps}>1-\delta\mu_{\theta}[-i\delta,i\delta]$.
Since $L_{p'}(\N)\lel [L_{q'}(\N),L_{r'}(\N)]_{\theta}$, we may find an
approximately norm-attaining element for $T_p(\psi^{\frac1p})$ as a simple
element of the interpolation space.  That is, we find an analytic function
$g:\mathfrak{S} \to L_{r'}(\N)$, continuous on the
boundary and vanishing at $\infty$ such that
  \[ \noo g(it)\bar{\phi}^{\frac{1}{r}-\frac{1}{q}}\rrm_{q'}\le
  (1+\eps), \qquad \noo g(1+it)\rrm_{r'} \kl (1+\eps),  \]
and
  \[ 1-\eps \le
  |tr_\N (g(\theta) \bar{\phi}^{\frac{1}{r}-\frac{1}{p}} T_p(\psi^{\frac1p}))|
  \pl . \]
Let $h$ be the function given by Lemma \ref{exKos}, and set
  \[ F(z)\lel tr_\N (g(z) \bfi^{\frac{1}{r}}\pi(h(z)))
  \lel tr_\N (g(z)T_r(\phi^{\frac{1}{r }} h(z))) \pl .\]
Since $z\mapsto \phi^{\frac1r} h(z)$ is analytic in the strip and
continuous at the boundary and bounded, we know that $F$ is
analytic, continuous at the boundary and vanishes at $\infty$.
Therefore, we deduce that
  \begin{align}
\label{estimate} 1-\eps &\le |tr(g(\theta) \bfi^{\frac1r - \frac1p}
T_p(\phi^{\frac1p} h(\theta))| = |tr(g(\theta)
\bfi^{\frac{1}{r}}\pi(h(\theta))| = |F(\theta)|\\
\notag &\lel \left|\intt F(z) d\mu_{\theta}(z) \right| \kl
   \intt_{i\rz} |F(it)| d\mu_{\theta}(it) +
   \intt_{1+i\rz} |F(1+it)| d\mu_{\theta}(it) \pl .
   \end{align}
For all $t$ we have
  \begin{align*}
   |F(it)| &=|tr_\N (g(it)T_r(\phi^{\frac1r} h(it) )| \lel
   |tr_\N (g(it) \bfi^{\frac1r}\pi(h(it)))| \lel
    |tr_\N ( g(it) \bfi^{\frac{1}{r}-\frac{1}{q}} T_q(\phi^{\frac1q} h(it) )|
    \\
   &\le \noo g(it)\bfi^{\frac{1}{r}-\frac{1}{q}} \rrm_{q'} \noo
\phi^{\frac{1}{q}}  h(it)\rrm_q \kl (1+\eps) \pl .
  \end{align*}
However, for $|t|\le \delta$, we note by \eqref{E:notisom} that the inequality
is stronger:
  \begin{align*}
   |F(it)| &=|tr_\N (g(it)T_r(\phi^{\frac1r} h(it) )| \lel
   |tr_\N (g(it) \bfi^{\frac1r}\pi(h(it)))| \lel
    |tr_\N ( g(it) \bfi^{\frac{1}{r}-\frac{1}{q}} T_q(\phi^{\frac1q} h(it) )|
    \\
   &\le \noo g(it)\bfi^{\frac{1}{r}-\frac{1}{q}} \rrm_{q'} \noo
T(\phi^{\frac{1}{q}}  h(it))\rrm_q \kl (1+\eps) \noo
T(\phi^{-\frac{it}{s}} \psi^{\frac{it}{s}} \psi^{\frac1q})\rrm_q \le (1 +
\eps)(1-\delta) \pl .
  \end{align*}
Similarly, we find
  \[ |F(1+it)|\kl 1+\eps \pl .\]
By \eqref{estimate}, this implies
  \[ {1-\eps} \kl (1 +
\eps)[(1-\delta)\mu_{\theta}([-i\delta,i\delta])+\mu_{\theta}(\partial
  \mathfrak{S} \setminus [-i\delta,i\delta])]
  \lel (1 + \eps)[1-\delta\mu_{\theta}([-i\delta,i\delta])] \pl ,\]
contradicting the choice of $\eps$.  This shows that $\noo T_q(\psi^{\frac1q})
\rrm_q=1$.  The argument for $\noo T_r(\psi^{\frac1r}) \rrm_q=1$
is similar.

It was shown in \cite{Ray} that the Mazur map $\psi\mapsto \psi^{\frac1q}$ is
continuous. So the above argument applies to a dense set of $L_q(M)_+$.  Since
$T_q$ is a contraction, we deduce for all positive elements $\psi^{\frac1q}$
  \[ \noo T_q(\psi^{\frac1q})\rrm_{L_q(\N)} \lel \|\psi^{\frac1q}\|_{L_q(\M)}
\pl .\]
For an arbitrary element $h \in L_q(\M)$, we write $h v^* = |\xi^*|$
where $v$ is a partial isometry. Then
  \[ \|h\|_q = \noo h^* \rrm_q \lel \noo \,|h^*|\,\rrm_q \lel \noo
T_q(|h^*|)\rrm_q \lel
  \noo T_q(h) \pi(v^*)\rrm_q \leq \noo T_q(h)\rrm_q \pl .\]
Therefore $T_q$ is an isometry. The same argument applies for
$T_r$. \qd

\begin{lemma}\label{41} Keep the notations of the previous proposition, and let $T_4$ be a $2$-isometry. Then for all $x\in \M$,
  \[ \noo \bfi^{\frac14}\pi(x)\bfi^{\frac14}\rrm_2 \lel
  \noo \phi^{\frac14} x \phi^{\frac14}\rrm_2 \pl .\]
\end{lemma}
\begin{proof} Let $x$ be a positive element, then we know that
  \[ \noo \bfi^{\frac14}\pi(x)\bfi^{\frac14}\rrm_2
  \lel \noo \bfi^{\frac14}\pi(x^{\frac12}) \rrm_4^2
  \lel \noo \phi^{\frac14} x^{\frac12} \rrm_4^2
  \lel \noo \phi^{\frac14}x\phi^{\frac14}\rrm_2 \pl .\]
Given an arbitrary element $x\in \M$ of norm less than one, the matrix
  \[ \tilde{x}\lel \left [ \begin{matrix} 1 & x\\ x^* & 1
  \end{matrix} \right ] \]
is positive. We apply the observation above to $\tilde{x}$ and
$\text{tr}_2\ten \phi$ as a normal faithful state on $M_2(\M)$ and deduce
  \[ \noo (\text{tr}_2\ten \bfi)^{\frac14}
\pi_2(\tilde{x})(\text{tr}_2\ten\bfi)^{\frac14}\rrm_2
  \lel \noo (\text{tr}_2\ten\phi)^{\frac14}  \tilde{x}
(\text{tr}_2\ten\phi)^{\frac14}\rrm_2
\pl .\]
Since $\noo\phi^{\frac 14} x^* \phi^{\frac 14}\rrm_2^2 =
\noo\phi^{\frac 14} x \phi^{\frac
14}\rrm_2^2$, we have
\begin{eqnarray*}
\noo (\text{tr}_2\ten \phi)^{\frac14}  \tilde{x}
(\text{tr}_2\ten\phi)^{\frac14}\rrm_2^2
&\lel& (\text{tr}_2 \otimes tr_\M) \left (
\left [\begin{matrix}\phi^{\frac 14} & 0 \\0 & \phi^{\frac
14}\end{matrix}\right ]
\left [\begin{matrix}1 & x^* \\x & 1\end{matrix}\right ]
\left [\begin{matrix}\phi^{\frac 12} & 0 \\0 & \phi^{\frac
12}\end{matrix}\right ]
\left [\begin{matrix}1 & x \\x^* & 1\end{matrix}\right ]
\left [\begin{matrix}\phi^{\frac 14} & 0 \\0 & \phi^{\frac
14}\end{matrix}\right ]
\right ) \\
&\lel& (\text{tr}_2 \otimes tr_\M) \left (
\left [\begin{matrix}\phi + \phi^{\frac 14} x^* \phi^{\frac 12} x
\phi^{\frac 14}&
\phi^{\frac 34}x \phi^{\frac 14}+ \phi^{\frac 14} x^* \phi^{\frac 34}\\
\phi^{\frac 34}x \phi^{\frac 14}+ \phi^{\frac 14} x^* \phi^{\frac 34} &
\phi + \phi^{\frac 14} x \phi^{\frac 12} x^* \phi^{\frac
14}\end{matrix}\right ]
\right ) \\
&\lel& {\frac 12} \left ( 2 tr_\M (\phi) +
tr_\M (\phi^{\frac 14} x^* \phi^{\frac 12} x \phi^{\frac 14})
+ tr_\M (\phi^{\frac 14} x \phi^{\frac 12} x^* \phi^{\frac 14}) \right ) \\
&\lel& 1 + \noo\phi^{\frac 14} x \phi^{\frac 14}\rrm_2^2.
\end{eqnarray*}
Since $\bfi(\pi(1))=1$, the same calculation shows that
\[  \noo (\text{tr}_2\ten \bfi)^{\frac14}  \pi_2(\tilde{x})
(\text{tr}_2\ten \bfi)^{\frac14}\rrm_2^2 \lel \left (1+\noo
\bfi^{\frac14}\pi(x)\bfi^{\frac14}\rrm_2^2 \right) \pl .\]
This proves the assertion. \qd

At this point we recall how complete isometries of $L_p$-spaces can be constructed from either *-isomorphisms or conditional expectations.  Details for the constructions of the next two paragraphs can be found in \cite{HRS}, \cite{JX}, and \cite{Sherman}.

Let $\pi: \M_1 \to \M_2$ be a surjective *-isomorphism, and let $\phi \in (\M_1)_*^+$ be faithful.  Then there is an associated completely isometric isomorphism $\pi_p: L_p(\M_1) \overset{\sim}{\to} L_p(\M_2)$, densely defined by
$$\phi^{\frac1p}x \mapsto (\phi \circ \pi^{-1})^{\frac1p}\pi(x), \qquad x \in \M_1.$$

Now consider the situation where $\M_1$ is a conditioned $\sigma$-finite subalgebra of $\M_2$, so that there are a normal *-isomorphism $\iota: \M_1 \hookrightarrow \M_2$ (thought of as the identity map) and a normal conditional expectation $E: \M_2 \twoheadrightarrow \M_1$.  By interpolation one can find a complete isometry $\iota_p: L_p(\M_1) \hookrightarrow L_p(\M_2)$ and a complete contraction $E_p: L_p(\M_2) \twoheadrightarrow \iota_p(L_p(\M_1))$.  Both $\iota_p$ and $E_p$ are completely positive.  With $\phi \in (\M_1)_*^+$ faithful, $\iota_p$ is densely defined by
$$\phi^{\frac1p} x \mapsto (\phi \circ E)^{\frac1p} x, \qquad x \in \M_1.$$

\begin{theorem}\label{T:sigmafinite} Let $1\le p\neq 2 <\infty$.  A linear map $T:L_p(\M)\to
L_p(\N)$ is a $2$-isometry if and only if there exist a normal $^*$-monomorphism
$\pi:\M\to \N$, a partial isometry $w \in \N$,
and a normal conditional expectation $E:\N\to \pi(\M)$ such that
$\pi(1) = w^*w$ and
\begin{equation} \label{E:phiform}
  T(\phi^{\frac1p}x)\lel w (\phi\circ \pi^{-1} \circ E)^{\frac1p}\pi(x), \qquad \forall x \in \M.
\end{equation}
Under these conditions, the range $T(L_p(\M))$ is completely contractively complemented, and $T$ is a module map:
$$T(h x) = T(h) \pi(x), \qquad \forall x \in \M, \: h \in L_p(\M).$$
\end{theorem}

\begin{proof} First note that maps of the form \eqref{E:phiform} are always complete isometries with completely contractively complemented range: they decompose as $\pi_p$, then $\iota_p$, then left multiplication by $w$.  Each of these is a complete isometry, and the image is the range of the complete contraction $L_p(\M_2) \ni h \mapsto w E_p(w^* h)$.

Now we turn to the derivation of \eqref{E:phiform} for 2-isometries.  The case $p=1$ has been proved  in Proposition \ref{T:lone}.
To prove the result for general $p$,  we may first apply Proposition \ref{T:mod} and assume that
  \[ T(\phi^{\frac1p}x)\lel \bfi^{\frac1p} \pi(x)\pl .\]
We can obtain the result for general $T$ by multiplying by an appropriate
partial isometry.

According to Lemma \ref{comp1}, we also know that $\phi\lel
\bar{\phi}\circ \pi$.
If $T$ is a $2$-isometry for some $1< p \neq 2 < \infty$, then we can
claim from
Corollary \ref {dual1} and Proposition \ref {extra} that $T_4$  is a
$2$-isometry.
Then  Lemma \ref{41} shows  that
  \[ \noo \phi^{\frac14}x\phi^{\frac14}\rrm_2 \lel
  \noo \bfi^{\frac14}\pi(x)\bfi^{\frac14}\rrm_2 \pl . \]
By polarization, we deduce
  \[ (\phi^{\frac14}x\phi^{\frac14},\phi^{\frac14}y\phi^{\frac14})
  \lel
  (\bfi^{\frac14}\pi(x)\bfi^{\frac14},\bfi^{\frac14}\pi(y)\bfi^{\frac14})\pl
  \]
for all $x,y\in \M$.  This means that
  \[ tr_\M(\phi^{\frac12}x\phi^{\frac12}y^*)\lel
   tr_\N(\bfi^{\frac12}\pi(x)\bfi^{\frac12}\pi(y)^*) \pl .\]
In other words the sesquilinear selfpolar form \cite{Woronowicz}
$s_{\bfi}(x,y) \triangleq tr(\bfi^{\frac12}x\bfi^{\frac12}y^*)$ satisfies
\begin{equation} \label{E:hs}
 s_{\bfi}|_{\pi(\M)\times \pi(\M)}\lel s_{\bfi|_{\pi(\M)}} \pl .
\end{equation}
Now by a result of Haagerup and St\o rmer \cite[Theorem 4.2]{Haa-St}, equation \eqref{E:hs} and the equality $\pi(1) = s(\bfi)$ imply the existence of a faithful normal conditional expectation $F:\pi(1)\N\pi(1) \to \pi(\M)$ such that $\bfi = \bfi \circ F$.  Defining $E: x \mapsto F(\pi(1)x \pi(1)), \: x \in \N,$ we have
$$\bfi = \bfi \circ E = \phi \circ \pi^{-1} \circ E,$$ which establishes \eqref{E:phiform}.
\qd

\begin{proof}[Proof of Theorem \ref{main}]
The main difference between Theorem \ref{main} and Theorem \ref{T:sigmafinite} is that $T$ is described on the spanning set of positive vectors $\{ \varphi^{\frac1p} \mid \varphi \in \M_*^+\}$ instead of the dense set $\{\phi^{\frac1p}x \mid x \in \M\}$ (and $T$ becomes everywhere-defined by linearity).  The equation \eqref{comp} is a noncommutative version of \eqref{E:banach}, so that $T$ may be naturally viewed as a ``noncommutative weighted composition operator".  To establish \eqref{comp}, it is sufficient to show that the data $\pi, E, w$ of Theorem \ref{T:sigmafinite} do not depend on the choice of $\phi \in \M_*^+.$

For this, choose an arbitrary faithful $\psi \in \M_*^+$ satisfying
$\psi^{\frac2p} \le C \phi^{\frac2p}$ for some $C<\infty$.  (Again, this set
is dense.)  The assumption means that $d \triangleq (D\phi:D\psi)_{i/p}$ exists in $\M$, and by analytic continuation we have $\phi^{\frac1p}d = \psi^{\frac1p}$.  (In fact the equation $d = \phi^{-\frac1p} \psi^{\frac1p}$ is justified rigorously in \cite{Sherman3}.)  Cocycles are functorial with
respect to normal *-isomorphisms, and they are invariant when weights are
precomposed with a normal conditional expectation \cite[Corollary
IX.4.22]{T2}, so
$$\pi(d) = (D(\phi \circ \pi^{-1} \circ E): D(\psi \circ \pi^{-1} \circ
E))_{i/p}.$$
Then we have
$$T(\psi^{\frac1p} ) = T(\phi^{\frac1p} d) = w (\phi \circ \pi^{-1} \circ E) ^{\frac1p} \pi(d) = w (\psi \circ \pi^{-1} \circ
E)^{\frac1p}.$$
By density, this establishes \eqref{comp}.  See also \cite[Section
6]{Sherman}.

Finally we note that the formulation of Theorem \ref{main} does not require the $\sigma$-finiteness of $\M$.  For each $q$ in the net of
$\sigma$-finite projections, the restricted isometry $T:L_p(q\M q) = qL_p(\M)q
\to L_p(\N)$ is of the form \eqref{comp} for some $w_q, \pi_q, E_q$.  One
only needs to glue them all together.  This can be done in exactly the same way
as in the proof of Theorem \ref{P.Jisometry}.

\end{proof}

\section{Concluding remarks}

\noindent \textbf{Remark 1.} If, in the setup of Theorem \ref{main}, we require that $ w^* w = \pi(1)$ is the support of $E$,
then $\pi$,  $w$, and $E$ are uniquely determined.  For suppose
that $\pi_1, E_1, w_1$ also define the 2-isometry $T$, and assume that $s(E_1) =
\pi_1(1)$.  Then for all $\varphi \in \M_*^+$,
\begin{align*}
w (\varphi \circ \pi^{-1} \circ E)^{\frac1p} &= w_1(\varphi \circ \pi_1^{-1} \circ
E_1)^{\frac1p} \Rightarrow \varphi \circ \pi^{-1} \circ E = \varphi \circ \pi_1^{-1}
\circ E_1 \\
&\Rightarrow \pi^{-1} \circ E = \pi_1^{-1} \circ E_1 \Rightarrow \text{id} =
\pi^{-1} \circ E \circ \pi_1 \Rightarrow \pi = E \circ \pi_1.
\end{align*}
Now take a projection $p \in \M$ and calculate
$$E(\pi(p)\pi_1(p)\pi(p)) = \pi(p)E(\pi_1(p))\pi(p) = \pi(p)^3 = \pi(p) =
E(\pi(p)).$$
Since $E$ is faithful on $\pi(1)\N\pi(1)$, we must have $\pi(p)\pi_1(p)\pi(p)
= \pi(p)$, or $\pi_1(p) \ge \pi(p)$.  Reversing the argument proves the
equality of $\pi$ and $\pi_1$, and the other data must be equal also (subject
to $w^*w = \pi(1) = w_1^*w_1$).

\bigskip

\noindent \textbf{Remark 2.} Given such a 2-isometry $T$ decomposed as in Theorem \ref{main}, the associated map $$S: L_p(\M) \to L_p(\N), \qquad h \mapsto w^*T(h),$$
is a completely positive complete isometry whose range is completely positively and completely contractively complemented.  If $T$ is positive, then $w = \pi(1)$ and $T=S$ (and in particular $T$ is already completely positive).

\bigskip

\noindent \textbf{Remark 3.} The arguments in section \ref{S:2isom} can be used to establish the following result, which seems to be of independent interest.  The main step strengthens Lemma \ref{41} by only requiring $T$ to be an isometry (i.e. not a 2-isometry).

\begin{prop}
Let $\M \subset \N$ be a unital inclusion of von Neumann algebras with $\phi \in \M_*^+, \bfi \in \N_*^+$, both faithful, satisfying $\bfi \mid_\M = \phi$.  Assume that the map
$$T_p: L_p(\M) \to L_p(\N); \qquad \phi^{\frac1p} x \mapsto \bfi^{\frac1p} x, \qquad x \in \M,$$
is isometric for some $1 \le p \ne 2 < \infty$.  Then there exists a faithful normal conditional expectation $E: \N \to \M$ such that $\bfi = \phi \circ E$.
\end{prop}

\begin{proof}
It follows from Lemma \ref{intpol}, Corollary \ref{dual1}, and Proposition \ref{extra} that if $T_p$ is isometric for one $p$ in the given range, it is isometric for all.  In particular, $T_4$ is isometric, so that for any $\M \ni y \ge 0$,
\begin{equation} \label{E:t4}
\|\phi^{\frac14} y \phi^{\frac14}\| = \|\phi^{\frac14} y^{\frac12} \|^2 = \|\bfi^{\frac14} y^{\frac12} \|^2 = \|\bfi^{\frac14} y \bfi^{\frac14}\|.
\end{equation}

Suppose that $y \in \M$ is only self-adjoint.  For all scalar $t \ge \|y\|$, \eqref{E:t4} implies
\begin{align*}
\|\phi^{\frac14} y \phi^{\frac14}\|^2 + 2t \phi(y) + t^2\|\phi\| &= \|\phi^{\frac14} (y + t1) \phi^{\frac14}\|^2\\
&= \|\bfi^{\frac14} (y + t1) \bfi^{\frac14}\|^2 \\
&= \|\bfi^{\frac14} y \bfi^{\frac14}\|^2 + 2t \bfi (y) + t^2\|\bfi\|.
\end{align*}
Since these are polynomials in $t$ which agree on a half-line, they must have the same constant terms.  We conclude that \eqref{E:t4} holds for all self-adjoint $y \in \M$.

In other words the map $\phi^{\frac14} y \phi^{\frac14} \mapsto \bfi^{\frac14} y \bfi^{\frac14}$ (with $y \in \M_{sa}$) is isometric between the real Hilbert spaces $L_2(\M)_{sa}$ and $L_2(\N)_{sa}$.  It therefore preserves inner products, which is the same as saying that the self-polar form $s_{\phi}$ agrees with $s_{\bfi}$ restricted to $\M \times \M$.  Then the Haagerup-St\o rmer result \cite[Theorem 4.2]{Haa-St} again implies the existence of a faithful normal conditional expectation $E: \N \to \M$ satisfying $\bfi = \phi \circ E$.
\end{proof}


\begin{thebibliography}{99}

\bibitem{Arazy} J. Arazy,
\emph{The isometries of $C_p$},
Israel J. Math. \textbf{22} (1975), 247--256.

\bibitem{Banach} S. Banach,
\emph{Th\'{e}orie des operations lin\'{e}aries},
Warsaw, 1932.

\bibitem{BL} J. Bergh and J. L\"{o}fstr\"{o}m,
\emph{Interpolation Spaces: an Introduction},
Springer-Verlag, Berlin-New York, 1976.

\bibitem{Broise} M. M. Broise,
\emph{Sur les isomorphismes de certaines alg\`{e}bres de von Neumann},
Ann. Sci. \'{E}cole Norm. Sup. (4) {\bf 83} (1966), 91--111.

\bibitem{Connesbook} A. Connes,
\emph{Noncommutative Geometry},
Harcourt Brace \& Co., San Diego, 1994.

\bibitem{Fidaleo} F. Fidaleo,
\emph{Canonical operator space structures on non-commutative $L\sp p$ spaces}, J. Funct. Anal. \textbf{169} (1999), no. 1, 226--250.

\bibitem{FJ} R. Fleming and J. Jamison,
\emph{Isometries on Banach Spaces: Function Spaces},
Chapman \& Hall/CRC Monographs and Surveys in Pure and Applied
Mathematics \textbf{129}, Boca Raton, 2003.

\bibitem{Haagerupwt} U. Haagerup,
\emph{Normal weights on W*-algebras},
J. Funct. Anal. \textbf{19} (1975), 302--317.

\bibitem{Haagerup} U. Haagerup,
``$L^p$-spaces associated with an arbitrary von Neumann algebra" in \emph{Alg\`{e}bres d'op\'{e}rateurs et leurs applications en
physique math\'{e}matique}, CNRS \textbf{15} (1979), 175--184.

\bibitem{Haa-St} U. Haagerup and E. St\o rmer,
\emph{Positive projections of von Neumann algebras onto JW-algebras},
Rep. Math. Phys. \textbf{36}, no. 2/3 (1995), 317--330.

\bibitem{HRS} U. Haagerup, H. Rosenthal, and F. Sukochev,
\emph{Banach embedding properties of non-commutative $L_p$-spaces},
Mem. Amer. Math. Soc. \textbf{163} (2003), no. 776.

\bibitem{HvN} P. Halmos and J. von Neumann,
\emph{Operator methods in classical mechanics II},
Ann. of Math \textbf{43} (1942), 332--350.

\bibitem{Ju} M. Junge,
\emph{Doob's inequality for non-commutative martingales},
J. Reine Angew. Math. \textbf{549} (2002), 149--190.

\bibitem{JRX} M. Junge,  Z-J. Ruan and Q. Xu,
\emph {Rigid ${\OL}_p$ structures of non-commutative $L_p$-spaces
associated with hyperfinite von Neumann algebras},
preprint.

\bibitem{JS} M. Junge and D. Sherman, 
\emph{Noncommutative $L^p$ modules},
J. Operator Theory, to appear.

\bibitem{JX} M. Junge and Q. Xu,
\emph{Noncommutative Burkholder/Rosenthal inequalities},
Ann. Probability, \textbf{31} (2003), 948--995.

\bibitem{Kat1} A. Katavolos,
\emph{Isometries of non-commutative $L^p$-spaces},
Canad. J. Math. \textbf{28} (1976), 1180--1186.

\bibitem{Kat2} A. Katavolos,
\emph{Are non-commutative $L^p$-spaces really non-commutative?}
Canad. J. Math. \textbf{33} (1981), 1319--1327.

\bibitem{Kat3} A. Katavolos,
\emph{Non-commutative $L^p$-spaces II},
Canad. J. Math. \textbf{34} (1982), 1208--1214.

\bibitem{Kirchberg} E. Kirchberg,
\emph{On nonsemisplit extensions, tensor products and exactness of group C*-algebras},
Invent. Math. \textbf{112} (1993), 449--489.

\bibitem{Ko1} H. Kosaki,
\emph{Applications of uniform convexity of non-commutative $L^p$-spaces},
Trans. Amer. Math. Soc. \textbf{283} (1984), 265--282.

\bibitem{Ko2} H. Kosaki,
\emph{Applications of the complex interpolation method to a von Neumann algebra: non-commutative $L^p$-spaces},
J. Funct. Anal \textbf{56} (1984), 29--78.

\bibitem{Lacey} H. Elton Lacey,
\emph{The Isometric Theory of Classical Banach Spaces},
Springer-Verlag, New York, 1974.

\bibitem{Lamperti} J. Lamperti,
\emph{On the isometries of certain function spaces},
Pacific J. Math. \textbf{8}(1958), 459--466.

\bibitem{Nelson} E. Nelson,
\emph{Notes on non-commutative integration},
J. Funct. Anal. \textbf{15} (1974), 103--116.

\bibitem{vNSchatten} J. von Neumann and R. Schatten,
\emph{The cross-space of linear transformations. III.}
Ann. of Math. (2) \textbf{49} (1948), 557--582.

\bibitem{NO} P. Ng and N. Ozawa,
\emph{A characterization of completely 1-complemented subspaces of noncommutative $L_1$-spaces},
Pacific J. Math. \textbf{205} (2002), no. 1, 171--195.

\bibitem{Ozawa} N. Ozawa,
\emph{Almost completely isometric embeddings between preduals of von Neumann algebras},
J. Funct. Anal. \textbf{186} (2001), no. 2, 329--341.

\bibitem{PT} G. K. Pedersen and M. Takesaki,
\emph{The Radon-Nikodym theorem for von Neumann algebras},
Acta Math. \textbf{130} (1973), 53--87.

\bibitem{PiLp} G. Pisier,
\emph{Non-commutative vector valued $L_p$-spaces and completely $p$-summing maps},
Ast\'{e}risque No. 247 (1998).

\bibitem{Pibook} G. Pisier,
\emph{An Introduction to the Theory of Operator Spaces},
London Mathematical Society Lecture Note Series \textbf{294}, Cambridge University Press, Cambridge, 2003.

\bibitem{PX} G. Pisier and Q. Xu,
``Non-commutative $L^p$ spaces" in \emph{Handbook of the Geometry of Banach Spaces, Volume 2}, North-Holland, New York, 2003.

\bibitem{Ray} Y. Raynaud,
\emph{On ultrapowers of non commutative $L_p$ spaces},
J. Operator Theory \textbf{48} (2002), 41--68.

\bibitem{RX} Y. Raynaud and Q. Xu,
\emph{On subspaces of noncommutative $L_p$ spaces},
preprint.

\bibitem{Ruan1} Z-J. Ruan,
\emph{Subspaces of $C^{*}$-algebras},
J. Funct. Anal. \textbf {76} (1988), 217--230.

\bibitem{Russo} B. Russo,
\emph{Isometries of $L^p$-spaces associated with finite von Neumann algebras},
Bull. Amer. Math. Soc. \textbf{74} (1968), 228--232.

\bibitem{Segal} I. Segal,
\emph{A non-commutative extension of abstract integration},
Ann. of Math. \textbf{57} (1953), 401--457.

\bibitem{Sherman} D. Sherman,
\emph{On the structure of isometries between noncommutative $L^p$ spaces}, preprint.

\bibitem{Sherman2} D. Sherman,
\emph{Noncommutative $L^p$ structure encodes exactly Jordan structure}, preprint.

\bibitem{Sherman3} D. Sherman,
\emph{Applications of modular algebras}, in preparation.

\bibitem{Stratila}  S. Str\v{a}til\v{a},
\emph{Modular Theory in Operator Algebras},
Abacus Press, Kent, 1981.

\bibitem{Ta} M. Takesaki,
\emph{Conditional expectations in von Neumann algebras},
J. Funct. Anal. \textbf{9} (1972), 306--321.

\bibitem{T1} M. Takesaki,
\emph{Theory of operator algebras I},
Springer-Verlag, New York, 1979.

\bibitem{T2} M. Takesaki,
\emph{Theory of operator algebras II},
Springer-Verlag, New York, 2002.

\bibitem{T3} M. Takesaki, 
\emph{Theory of operator algebras III},
Springer-Verlag, New York, 2002.

\bibitem{Tam} P. K. Tam,
\emph{Isometries of $L^p$-spaces associated with semifinite
von Neumann algebras},
Trans. Amer. Math. Soc. \textbf{154} (1979), 339--354.

\bibitem{Terp1} M. Terp,
\emph{$L^p$ spaces associated with von Neumann algebras},
notes, Copenhagen Univ., 1981.

\bibitem{Terp2} M. Terp,
\emph{Interpolation spaces between a von Neumann algebra and its predual},
J. Operator Theory, \textbf{8} (1982), 327--360.

\bibitem{W1} K. Watanabe,
\emph{Problems on isometries of non-commutative $L^p$-spaces},
Contemp. Math. \textbf{232} (1989), 349--356.

\bibitem{W2} K. Watanabe,
\emph{On the structure of non-commutative $L^p$-isometries}, preprint.

\bibitem{Woronowicz} S. L. Woronowicz,
\emph{Selfpolar forms and their applications to the C*-algebra theory},
Rep. Mathematical Phys. \textbf{6} (1974), no. 3, 487--495.

\bibitem{Yamagami} S. Yamagami,
\emph{Algebraic aspects in modular theory},
Publ. RIMS \textbf{28} (1992), 1075--1106.

\bibitem{Yeadon} F. J. Yeadon,
\emph{Isometries of non-commutative $L^p$-spaces},
Math. Proc. Camb. Phil. Soc. \textbf{90} (1981), 41--50.

\end{thebibliography}
\end{document}